\documentclass[12pt]{amsart}
\usepackage{amsmath}
\usepackage{amssymb}
\usepackage{amsthm}
\usepackage{mathrsfs}
\usepackage{comment}
\usepackage{hyperref}

\usepackage[all,cmtip]{xy}\usepackage{xcolor}
\usepackage{enumerate}
\usepackage{bm}
\usepackage{dsfont}
\usepackage{mathtools}
\usepackage[cal=euler]{mathalfa}
\usepackage[top=1in, bottom=1.25in, left=1.25in, right=1.25in]{geometry}
\usepackage{parskip}

\hypersetup{colorlinks=true,linkcolor=magenta,citecolor=blue}

\DeclareMathAlphabet{\mathpzc}{OT1}{pzc}{m}{it}

\theoremstyle{plain}
\newtheorem{theorem}{Theorem}[section]
\newtheorem*{theorem*}{Theorem}

\newtheorem{lemma}[theorem]{Lemma}
\newtheorem*{claim*}{Claim}
\newtheorem{proposition}[theorem]{Proposition}

\theoremstyle{definition}
\newtheorem{definition}[theorem]{Definition}

\newtheorem{example}[theorem]{Example}

\newtheorem{remark}[theorem]{Remark}

\numberwithin{equation}{section}
\numberwithin{figure}{section}

\newcommand{\normKF}{N_{K/F}}
\newcommand{\normKFim}{N_{K/F}(K^\times)}

\newcommand{\N}{\mathbb{N}}
\newcommand{\F}{\mathbb{F}}

\newcommand{\ep}{\varepsilon}

\newcommand{\ignore}[1]{}

\author{S. Pumpl\"un}
\email{susanne.pumpluen@nottingham.ac.uk}
\address{School of Mathematical Sciences\\
University of Nottingham\\
University Park\\
Nottingham NG7 2RD\\
United Kingdom
}

\keywords{$(f,\sigma)$-codes, $\sigma$-constacyclic codes, nonassociative algebra, coset coding, maximum rank distance codes, Construction A, lattice,  learning with errors, post-quantum cryptography, NCLWE}

\subjclass[2020]{Primary: 17A35; Secondary: 11T71, 94B40, 94B05}

\begin{document}

\title[Lattice codes]
{Using cyclic $(f,\sigma)$-codes over finite chain rings to construct $\mathbb{Z}_p$- and $\F_q[\![t]\!]$-lattices}

\begin{abstract}
We construct $\mathbb{Z}_p$-lattices and
$\F_q[\![t]\!]$-lattices from cyclic $(f,\sigma)$-codes over finite chain rings, employing quotients of natural nonassociative orders and  principal left ideals in carefully chosen nonassociative  algebras.
This approach generalizes the classical Construction A that obtains $\mathbb{Z}$-lattices from linear codes over finite fields or commutative rings to the nonassociative setting.  We mostly use proper nonassociative cyclic algebras that are defined over field extensions of $p$-adic fields. This means we focus on $\sigma$-constacyclic codes to obtain $\mathbb{Z}_p$-lattices, hence $\mathbb{Z}_p$-lattice codes.
We construct linear maximum rank distance (MRD) codes that are $\mathbb{Z}_p$-lattice codes employing  the left multiplication of a nonassociative algebra over a finite chain ring.

Possible applications of our constructions include post-quantum cryptography involving $p$-adic lattices, e.g. learning with errors, building rank-metric codes like MRD-codes, or $p$-adic coset coding, in particular wire-tap coding.

\end{abstract}

\maketitle

\section*{Introduction}

Since $p$-adic numbers and $p$-adic lattices are increasingly of interest in some areas of coding theory and post-quantum cryptography, e.g. for
$p$-adic lattice public-key encryption, cryptosystems, and signature schemes \cite{Deng, Detal, Z, ZDW, Z2024}, we  will explore some  lattice constructions over local nonarchimedean fields using nonassociative algebras over local nonarchimedean fields, their natural orders and their ideals, to generate new lattices.

The classical Construction A lifts a linear code over the finite field $\mathbb{F}_2$ to obtain a $\mathbb{Z}$-lattice \cite{CS}.
For suitable skew polynomials $f\in K[t;\sigma,\delta]$, automorphisms $\sigma$, and left $\sigma$-derivations $\delta$,
canonical generalizations of this classical constructions lift a cyclic  $(\overline{f},\overline{\sigma},\overline{\delta})$-code $\mathcal{C}_0$ over a finite field or more generally, a Galois ring, to obtain a $\mathbb{Z}$-lattice code that is a cyclic $(f,\sigma,\delta)$-code  \cite{DO, OS, P18.2}.
These generalizations employ both associative and nonassociative algebras  over number fields, their rings of integers,  their natural orders and ideals.

 In the following, let $K/F$ be a finite Galois field extension of local nonarchimedean fields.
Let $\mathcal{O}_F$ and $\mathcal{O}_K$ be the valuation rings (i.e., the rings of integers) of $F$, respectively  $K$, with residue fields $\overline{F}=\mathcal{O}_F/\mathfrak{p}$ and $\overline{K}=\mathcal{O}_K/\mathfrak{p}_K$, where $\mathfrak{p}_K$ denotes the unique maximal ideal of $\mathcal{O}_K$ and
$\mathfrak{p}$  the unique maximal ideal of $\mathcal{O}_F$.

In this paper we work with certain associative and nonassociative algebras over local nonarchimedean fields
 and lift cyclic $(\overline{f},\overline{\sigma}, \overline{\delta})$-codes $\mathcal{C}$  that are defined over finite chain rings $\mathcal{O}_K/\mathfrak{p}\mathcal{O}_K$.

There are two kinds of local nonarchimedean fields:  the finite algebraic extensions of the $p$-adic numbers $\mathbb{Q}_p$ (the \emph{$p$-adic fields}), and the fields $\F_q(\!(t)\!)$ of Laurent series over a finite field of order $q=p^j$ for some $j\geq 1$. 
Each of them yields a very different setup for possible code and lattice constructions.

We construct $\mathbb{Z}_p$-lattices, respectively $\F_q[\![t]\!]$-lattices, out of cyclic $(\overline{f},\overline{\sigma},\overline{\delta})$-codes $\mathcal{C}$ over the finite ring $S=\mathcal{O}_K/\mathfrak{p}\mathcal{O}_K$. We mostly focus on $\overline{\sigma}$-constacyclic codes $\mathcal{C}$ over $S$ when $K$ is a $p$-adic field, and use these to construct $\mathbb{Z}_p$-lattice codes in $\mathbb{Q}_p^N$.

For the special case of $\overline{\sigma}$-constacyclic codes, the idea of our approach can be summarized as follows: we first choose a finite cyclic Galois field extension $K/F$ of local nonarchimedean fields with Galois group generated by $\sigma$. We then choose a monic  $f\in \mathcal{O}_K[t;\sigma]$ which is irreducible in $K[t;\sigma]$, such that the polynomial we obtain by reading its coefficients modulo the ideal $\mathfrak{p}\mathcal{O}_K$ is hopefully reducible (as the monic right divisors of this polynomial correspond to the codes we use, so the more monic right divisors there are, the more codes we can construct). We  take a natural order $\Lambda$ in the nonassociative ring $K[t;\sigma]/K[t;\sigma]f$ which may depend on a choice of a maximal subfield when the ring is associative, and observe that the ring $\Lambda/\mathfrak{p}\Lambda$ is isomorphic to  the nonassociative  Petit algebra $$\overline{A}=(\mathcal{O}_K/\mathfrak{p}\mathcal{O}_K)[t;\overline{\sigma}]/(\mathcal{O}_K/\mathfrak{p}\mathcal{O}_K)[t;\overline{\sigma}]f.$$
 We then choose a principal left ideal $\mathcal{I}$  of $\Lambda$ generated by a monic polynomial, such that
$\mathfrak{p}\subset \mathcal{I}$.   For $g(t)=\sum_{i=0}^{m}a_i t^i \in K[t;\sigma]$, we define $\overline{g}(t)=\sum_{i=0}^{m}\overline{a_i}t^i\in (\mathcal{O}_K/\mathfrak{p}\mathcal{O}_K)[t;\overline{\sigma}]$
  with $\overline{a_i}= a_i+\mathfrak{p}\mathcal{O}_K$.  Then the image  $\Psi(\mathcal{I}/\mathfrak{p}\Lambda)$ under the canonical map
 $$\Psi:\Lambda\rightarrow   \overline{A}= (\mathcal{O}_K/\mathfrak{p}\mathcal{O}_K)[t;\overline{\sigma}]/(\mathcal{O}_K/\mathfrak{p}\mathcal{O}_K)[t;\overline{\sigma}] \overline{f},\quad g\mapsto \overline{g}$$
 is also a principal left ideal in the algebra $\overline{A}$,  generated by the monic polynomial  $\overline{f}$. It therefore
  corresponds to a cyclic
$(\overline{f},\overline{\sigma})$-code $\mathcal{C}$ over $S=\mathcal{O}_K/\mathfrak{p}\mathcal{O}_K$. Let
$\rho: \Lambda\rightarrow  \Lambda/\mathfrak{p}\Lambda \rightarrow  \Psi(\Lambda/\mathfrak{p}\Lambda)$
be the canonical projection $\Lambda\rightarrow  \Lambda/\mathfrak{p}\Lambda$ composed with $\Psi$. Depending on our choice  of the local nonarchimedean field $K$, the pre-image
$\rho^{-1}(\mathcal{C})$ 
is either a $\mathbb{Z}_p$-lattice, respectively, an $\F_q[\![t]\!]$-lattice.

After collecting  the lesser known results and definitions we need in Section \ref{sec:prel}  (among them the definition of a nonassociative Petit algebra $S_f$ using skew polynomials $f\in K[t;\sigma,\delta]$), we set up  Construction A  in Section \ref{sec:naturalI}. We explain how to obtain $\mathbb{Z}_p$-lattice codes, respectively,  $\F_q[\![t]\!]$-lattice codes, from cyclic $(\overline{f},\overline{\sigma})$-codes over $\mathcal{O}_K/\mathfrak{p}\mathcal{O}_K$ using a generalized Construction A in Section
\ref{sec:codes}.
These results are then generalized  in Section \ref{sec:naturalII} using more general Petit algebras $S_f$ constructed using
 monic irreducible  $f\in D[t;\sigma,\delta]$, where $D=(K/F,\rho,c)$ is a cyclic division algebra.

In Section \ref{sec:MRD}, we lift linear matrix codes to maximum rank distance  (MRD) codes that are lattice codes by computing the matrix representing the left multiplication in a division algebra $\overline{A}$. The
 linear code  $\mathcal{C}$ associated to a principal left ideal $\Psi(\mathcal{I}/\mathfrak{p}\Lambda)$ then defines a linear matrix code  via $\rho^{-1}(\mathcal{C})$ 
 This is an  MRD code over $\mathcal{O}_K$, which is a $\mathbb{Z}_p$-lattice, respectively, a $\F_q[\![t]\!]$-lattice.
 We sketch potential applications to learning with errors (LWE) in Section
\ref{sec:outlook}.

This paper  generalizes an  approach presented in \cite{DO, P18.2}. Its weakness is a lack of explicit examples for the lattice constructions we propose, which will be given in a subsequent  paper; here we focus on developing the theory. Examples will exploit recent results in the classification on nonassociative cyclic algebras over local nonarchimedean number fields \cite{NP}
 The idea of lattice encoding can be also applied for coset coding and wiretap coding, analogously as described in \cite[Sections 5.2, 5.3]{DO}.
We refer the reader to \cite{P18.2} for how to apply our results to space-time block coding.

Note that codes over finite chain rings, especially Galois rings, have been studied for instance in \cite{B, Cao, LL, MW}.

\section{Preliminaries}\label{sec:prel}

\subsection{Nonassociative algebras}

Let $R$ be a unital commutative ring and let $A\not=0$ be an
$R$-module. Then $A$ is a  \emph{(nonassociative) algebra} over $R$, if there exists an
$R$-bilinear map $A\times A\to A$, $(x,y) \mapsto x \cdot y$, usually denoted by juxtaposition $xy$,
which defines the nonassociative multiplication of $A$. An algebra is called a \emph{proper nonassociative algebra}, if we want to stress that $A$ is not associative. An algebra $A$ is called \emph{unital} if there is an element in $A$, denoted by 1, such that $1x=x1=x$ for all $x\in A$.
All algebras in this paper are unital algebras.

Define  $[x, y, z] = (xy) z - x (yz)$ for $x,y,z\in A$.
For an $R$-algebra $A$, the {\it left nucleus}, {\it middle nucleus} and {\it right nucleus} of $A$ are defined as ${\rm Nuc}_l(A) = \{ x \in A \, \vert \, [x, A, A]  = 0 \}$, ${\rm Nuc}_m(A) = \{ x \in A \, \vert \, [A, x, A]  = 0 \}$ and   ${\rm Nuc}_r(A) = \{ x \in A \, \vert \, [A,A, x]  = 0 \}$. The  {\it nucleus} of $A$ is given by
 ${\rm Nuc}(A) = \{ x \in A \, \vert \, [x, A, A] = [A, x, A] = [A,A, x] = 0 \}$ and
the {\it center} of $A$ is ${\rm C}(A)=\{x\in A\,|\, x\in \text{Nuc}(A) \text{ and }xy=yx \text{ for all }y\in A\}$
 \cite{Sch}.

   An algebra $A$ over a field $F$ is called a {\it division algebra}, if for any $a\in A$, $a\not=0$,
the left multiplication  with $a$, $L_a(x)=ax$,  and the right multiplication with $a$, $R_a(x)=xa$, are bijective.
Any division algebra is simple, that means has only trivial two-sided ideals.
A finite-dimensional algebra
$A$ is a division algebra over $F$ if and only if $A$ has no zero divisors.

If $R$ is a Noetherian integral domain with quotient field $F$ and $A$ a finite-dimensional unital $F$-algebra, then
 an $R$-\emph{lattice} in $A$ is an $R$-submodule $\Gamma$ of $A$ which is finitely generated and contains an $F$-basis of $A$.
  A multiplicatively closed $R$-lattice containing $1_A$
(this multiplication may be not associative) is called an $R$-\emph{order} in $A$.

\subsection{Skew polynomial rings}

Let $S$ be a unital associative ring, $\sigma:S\rightarrow S$ an injective ring homomorphism  and $\delta:S\rightarrow S$ a \emph{left $\sigma$-derivation}, i.e. $\delta$ is an additive map such that $\delta(ab)=\sigma(a)\delta(b)+\delta(a)b$
for all $a,b\in S$, implying $\delta(1)=0$. Put ${\rm Const}(\delta)=\{a\in S\,|\, \delta(a)=0\}$ and
${\rm Fix}(\sigma)=\{a\in S\,|\, \sigma(a)=a\}$.

 The \emph{skew polynomial ring} $R=S[t;\sigma,\delta]$
is the set of \emph{skew polynomials} $a_0+a_1t$ $+\dots +a_nt^n$
with $a_i\in S$, where addition is defined term-wise as we know it for classical polynomials,  and a skewed multiplication is defined via the rule
$ta=\sigma(a)t+\delta(a)$ for all $a\in S$
 (\cite{O1}, for properties see \cite{C, G6}).

 For $f(t)=a_0+a_1t+\dots +a_nt^n\in R$, $a_n\not=0$, define ${\rm deg}(f)=n$ and for $f=0$ define  ${\rm deg}(f)=-\infty$.
Then ${\rm deg}(fg)\leq{\rm deg} (f)+{\rm deg}(g)$
with equality when $f$ has an invertible leading coefficient,
 $g$ has an invertible leading coefficient,
and $S$ is a domain.

 An element $f\in R$ is \emph{irreducible} in $R$ if it is not a unit and it has no proper factors, i.e if there do not exist $g,h\in R$ with
 ${\rm deg}(g),{\rm deg} (h)<{\rm deg}(f)$ such that $f=gh$.

\subsection{Petit's algebras}\label{sec:S_f}

  Assume that the skew polynomial
  $$f(t)=\sum_{i=0}^{m}d_it^i\in R=S[t;\sigma,\delta]$$
   has an invertible leading coefficient $d_m$. Then
for all $g(t)\in R$ of degree $l> m$,  there exist  uniquely determined $r(t),q(t)\in R$ with
 ${\rm deg}(r)<{\rm deg}(f)$, such that
$g(t)=q(t)f(t)+r(t)$, e.g. \cite[Proposition 1]{P15}.

Let  $R_m=\{g\in R\,|\, {\rm deg}(g)<m\}$ and let ${\rm mod}_r f$ denote the remainder of right division by  $f$.
For $g,h\in R_m$  we define a multiplication via
$g\circ h=  gh \,\,{\rm mod}_r f$, where the right hand side is the remainder we obtain after dividing the product $gh$ in $R$ on the right by $f$. Then $(R_m,\circ)$
is a  unital nonassociative ring that we denote by $S_f$ or by $R/Rf$  (\cite{P15}, or \cite{P66} when $S$ is a division ring).
The ring $S_f$ is associative if and only if  $Rf$ is a two-sided ideal in $R$ (\cite[Theorem 4 (ii)]{P15}, or \cite[(1)]{P66} when $S$ is a division ring).

Define $S_0=\{a\in S\,|\, ah=ha \text{ for all } h\in S_f\}$, then $S_0$  is a commutative subring of $S$, and  $S_f$ is a unital nonassociative algebra  over $S_0$. We call $S_f$ a \emph{Petit algebra} as the construction goes back to Petit \cite{P66}.

 When $S_0$ is a field, $S$ is a division algebra,  and $S_f$ is a finite-dimensional vector space over $S_0$,
  then $S_f$ is a division algebra if and only if $f$ is irreducible in $R$ \cite[(9)]{P66}.

 If $S_f$ is a proper nonassociative algebra then $S\subset {\rm Nuc}_l(S_f)$ and $S\subset{\rm Nuc}_m(S_f)$,
${\rm Nuc}_r(S_f)=\{g\in R_m\,|\, fg\in Rf\}$ and $S_0$ is the center of $S_f$ \cite{P15}.
It is easy to see that
$$C(S)\cap {\rm Fix}(\sigma)\cap {\rm Const}(\delta)\subset S_0.$$
If $f$ is irreducible, $\delta=0$   and $S$ is a division ring then
$S_0={\rm Fix}(\sigma)\cap C(S)$ is the center of $S_f$ \cite[Theorem 8 (ii)]{P15}.

We call the  $S$-basis $1, t, t^2, \ldots, t^{n-1}$   the \emph{canonical basis} for the left $S$-module $S_f$.
The left multiplication with $0\not=a\in S_f$ in $S_f$, $L_a:S_f\rightarrow  S_f$, is an $S$-module endomorphism. We can hence express $L_a$ in matrix form with respect to the canonical basis of $S_f$ and obtain a map $\gamma: S_f \to {\rm End}_K(S_f), a\mapsto L_a$ which in turn
induces an injective $S$-linear map
$$\gamma: S_f \to {\rm Mat}_m(S), a\mapsto L_a \mapsto M(a).$$
This map is used when designing maximum rank distance or space-time block  codes which often employ  the following type of algebras which will also be the focus of this paper:

\begin{definition}
 Let $S/S_0$ be an extension of commutative unital rings and $G=\langle \sigma\rangle$ a finite cyclic group of order
$m$ acting on $S$ such that $S_0={\rm Fix}(\sigma)$. For any $c\in S$,
$$S_f=S[t;\sigma]/S[t;\sigma] (t^m-c)$$
is called a \emph{nonassociative cyclic algebra}  $(S/S_0,\sigma,c)$ \emph{of degree $m$} over $S_0$.
\end{definition}

 When $c \in S \setminus S_0$, then $(S/S_0,\sigma,c)$ has left and middle nucleus containing
$S$  and center $S_0$. When $S/S_0$ is a cyclic Galois field extension with Galois group generated by $\sigma$, and $c\in S_0^\times$, then $S[t;\sigma]/S[t;\sigma] (t^m-c)$ is a classical associative cyclic algebra.  When $c=0$, the algebra $S[t;\sigma]/S[t;\sigma] (t^m)$ is  commutative and associative and is the direct product of $m$ copies of $S$.

   If $S/S_0$ is a cyclic
 Galois field extension of degree $m$ with Galois group $\langle \sigma\rangle$ and $c\in S\setminus S_0$,
 then ${\rm Nuc}((S/S_0,\sigma,c))=S$. Moreover, for non-prime $m$, $(S/S_0,\sigma,c)$ is a division
 algebra for all choices of $c$ such that $1,c,\dots,c^{m-1}$ are linearly independent \cite{S12}, in particular when $m$ is prime.

\subsection{Cyclic $(f,\sigma,\delta)$-codes} \label{subsec:codes}

Let $S$ be a commutative unital ring. A \emph{linear code of length $m$ over $S$} is a submodule of the $S$-module $S^m$.

 Let $f\in S[t;\sigma,\delta]$ be monic of degree $m$ and $\sigma$ injective. We can associate to an element
$a(t)=\sum_{i=0}^{m-1}a_it^i$ in $S_f$ the  vector $(a_0,\dots,a_{m-1})\in S^m$.
 Conversely, for any linear code $\mathcal{C}$ of length $m$ 
 we have a set of skew polynomials
 $a(t)=\sum_{i=0}^{m-1}a_it^i\in S_f$ associated to the codewords $(a_0,\dots,a_{m-1})\in \mathcal{C}$.

A \emph{cyclic $(f,\sigma,\delta)$-code} $\mathcal{C}\subset S^m$ is a set consisting of the vectors
$(a_0,\dots,a_{m-1})$ obtained from all the elements $h=\sum_{i=0}^{m-1}a_it^i$
in a left principal ideal $S_f g$ where $S_f=S[t;\sigma,\delta]g/S[t;\sigma,\delta]f$, and $g$ is a monic
right divisor of $f$.
A code $\mathcal{C}$ over $S$ is called
   $\sigma$-\emph{constacyclic} if there is a non-zero $c\in S$ such that
   $$(a_0,\dots,a_{m-1})\in  \mathcal{C}\Rightarrow (\sigma(a_{m-1})c,\sigma(a_0),\dots,\sigma(a_{m-2}))\in  \mathcal{C}$$
   (note that $\sigma$-constacyclic codes are obtained when $f(t)=t^n-c\in S[t;\sigma]$).

 Every  right divisor $g$ of $f$ of degree $<m$ with an invertible leading coefficient generates a principal left ideal in $S_f$.
 Indeed, all the left ideals in $S_f$ which contain a non-zero  polynomial $g$
of minimal degree with invertible leading coefficient are principal left ideals, and $g$ is a right divisor of $f$ in $R$.  When $S$ is a division ring then all left ideals in $S_f$ are generated by some monic right divisor $g$ of $f$ (cf. \cite[Proposition 7]{P15}).
  Each principal left ideal generated by a monic right divisor of $f$ is
an $S$-module which is isomorphic to a submodule of $S^m$ and forms a code of length $m$ and dimension $m-{\rm deg}(g)$ \cite[Theorem  1]{BL13}.

For a field $K$, every skew polynomial ring $K[t;\sigma,\delta]$ can be made into either a twisted or a
differential polynomial ring by a linear change of variables  \cite[1.1.21]{J96}. When constructing linear codes, however,
it makes sense to study general skew polynomial rings to not limit the range of linear codes we are able to construct, for instance cf. \cite{BU14} for examples of this phenomenon.

\subsection{Local nonarchimedean fields} (cf. \cite[Chapter II]{Lang1994}; for examples cf. \cite{LMFDB}).

Let $F$ be a local nonarchimedean field with ring of integers $\mathcal{O}_F$ and let $\mathfrak{p}$ be the unique maximal ideal of $\mathcal{O}_F$.  The residue field $\overline{F}=\mathcal{O}_F/\mathfrak{p}$ of $F$ is a finite field of order $q=p^j$ for some $j\in \N$ and prime $p$ ($p$ is called the \emph{residual characteristic} of $F$).

\begin{example} If $F=\mathbb{F}_q(\!(t)\!)$ then $\mathcal{O}_F= \F_q[\![t]\!]$, $\mathfrak{p}=( t )$, and $\overline{ F }= \mathbb{F}_q$. If $F=\mathbb{Q}_p$ then
$$\mathbb{Q}_p^\times=\{ a=\sum_{i=N}^\infty a_i p^i : N\in \mathbb{N}, a_i \in \{0,1,\ldots,p-1\}\}
$$
with addition and multiplication computed ``with carrying" as one would for finite sums. Here the least $i$ such that $a_i\neq 0$ is the $p$-adic valuation of $a$, and in this case $\vert a\vert = p^{-i}$.  The set of all elements of valuation at least $0$ are the \emph{$p$-adic integers} $\mathcal{O}_F=\mathbb{Z}_p$.
\end{example}

\begin{lemma}[Hensel's Lemma] \label{le:Hensel} \cite[Thm 7.3]{Eisenbud1995}
Let $f(x)\in \mathcal{O}_F [x]$ be a polynomial and $f'(x)$ its formal derivative.  If $a\in \mathcal{O}_F $ satisfies
$$
f(a) \equiv 0 \mod f'(a)^2 \mathfrak{p},
$$
then there exists a unique $b\in \mathcal{O}_F $ such that
$$
f(b)=0  \text{ and } b\equiv a \mod f'(a) \mathfrak{p}.
$$
\end{lemma}

Finite algebraic extensions $F$ are again local nonarchimedean fields and they come in two forms:  There is a unique (Galois) \emph{unramified extension} for each degree $n$ which  we will denote by $K_n$, which is the splitting field of $x^{q^n}-x$.  We have $\overline{K_n}\cong \F_{q^n}$ and the maximal ideal $ \mathfrak{p}_{K_n}$ is generated over $\mathcal{O}_{L_n}$ by $\varpi$. The Galois group of $K_n/F$ is cyclic.

At the other extreme are the \emph{totally ramified extensions} $K/F$ of degree $n$.
 Their residue fields satisfy $\overline{ K} \cong \overline{ F}$ and $\varpi$ generates $ \mathfrak{p}_K^n$. Not all of them are Galois extensions.

A general finite algebraic extension $K$ of $F$ can be uniquely factored as an unramified extension $L/F$ of degree $f$ followed by a totally ramified extension $K/L$ of degree $e$; $e$ is called the \emph{ramification index} of $K/F$.  We say $K/F$ is \emph{tamely ramified} if $(e,p)=1$ and \emph{wildly ramified} if not.  The totally and tamely ramified extensions of degree $e$ are all obtained by adjoining to $F$ a root of $x^n-u\varpi$, for some choice of $u\in \mathcal{O}_F^\times/(\mathcal{O}_F^\times)^e$; these are Galois only if $F$ contains a primitive $e$th root of unity.

The Galois extensions $K/F$ of degree $n$  are in one-to-one correspondence with the subgroups of index $n$ of $F^\times$, via the map that assigns $K/F$ to the subgroup $\normKFim$, which is the image of the norm map \cite[XI \S4]{Lang1994}, hence $F^\times/\normKFim$ is a finite group of order $n$. 

If $K=K_n$ is the unramified extension of $F$, then  restricting the norm map to the ring of integers we obtain a surjective map
$\normKF : \mathcal{O}_{K}^\times \to \mathcal{O}_F^\times$,
 and the group $F^\times/\normKFim$ is represented by the elements
$$
\{ 1, \varpi, \varpi^2, \cdots, \varpi^{n-1} \},
$$
 while for a totally and tamely ramified extension $K/F$, $F^\times/\normKFim \cong \mathcal{O}_F^\times/(\mathcal{O}_F^\times)^n$  and Hensel's Lemma identifies the group $F^\times/\normKFim$ with $\overline{ F}^\times /(\overline{F}^\times)^n$.

In the following, we use the approach from  \cite[Section 2]{P18.2}, but employ local nonarchimedean fields instead of number fields.

\section{Natural orders in $S_f$ and their quotients by a prime ideal} \label{sec:naturalI}

\subsection{Quotients of rings on integers} \label{subsec:assumptions}

Let $K/F$ be a Galois extension of local nonarchimedean fields of degree $n$ and  $G={\rm Gal}(K/F)$.
Let $\mathcal{O}_F$ and $\mathcal{O}_K$ be the ring of integers of $F$, respectively  $K$.
Let $\mathfrak{p}_K$ be the unique maximal ideal of $\mathcal{O}_K$, $\overline{K}=\mathcal{O}_K/\mathfrak{p}_K$, and
$\mathfrak{p}$ be the unique maximal ideal of $\mathcal{O}_F$, so that $\overline{F}=\mathcal{O}_F/\mathfrak{p}=\mathbb{F}_{p^j}$  for some $j\in \N$ and prime $p$.

 Let $\pi:\mathcal{O}_K\rightarrow  \mathcal{O}_K/\mathfrak{p}\mathcal{O}_K$ be
the canonical projection. For every $\sigma\in G$, we have $\sigma(\mathfrak{p}\mathcal{O}_K)\subset \mathfrak{p}\mathcal{O}_K$ since
$\sigma|_F=id$. Thus $\sigma$ induces a ring homomorphism
$$\overline{\sigma}:\mathcal{O}_K/\mathfrak{p}\mathcal{O}_K \rightarrow  \mathcal{O}_K/\mathfrak{p}\mathcal{O}_K,
\quad a+\mathfrak{p}\mathcal{O}_K \mapsto \sigma(a)+\mathfrak{p}\mathcal{O}_K $$
satisfying $\pi\circ \sigma|_{\mathcal{O}_K}=\overline{\sigma} \circ \pi$ and ${\rm Fix}(\overline{\sigma})=\mathbb{F}_{p^j}$.
Suppose that $\delta$ is an $F$-linear left $\sigma$-derivation on $K$ such that
$\delta(\mathcal{O}_K)\subset \mathcal{O}_K$.
 Then $\delta$ induces a left $\overline{\sigma}$-derivation
 $$\overline{\delta}:\mathcal{O}_K/\mathfrak{p}\mathcal{O}_K \rightarrow  \mathcal{O}_K/\mathfrak{p}\mathcal{O}_K.$$

A general finite algebraic extension $K/F$ can be uniquely written as a tower of field extensions $F\subset E\subset K$ where $K/E$ is totally ramified of degree $e\geq 1$,  the \emph{ramification index}, and $E/F$ is unramified of degree $f\geq 1$.
We say $K/F$ is \emph{tamely ramified} if $(e,p)=1$ and \emph{wildly ramified} if not.  The totally and tamely ramified extensions of degree $e$ are all obtained by adjoining to $F$ a root of $x^n-u\varpi$, for some choice of $u\in (\mathcal{O}_F^\times/(\mathcal{O}_F^\times)^e$; these are Galois only if $F$ contains a primitive $e$th root of unity.

For $ \mathfrak{p}=(\varpi)$ we have $ \mathfrak{p}\mathcal{O}_E=\mathfrak{p}_E =(\varpi)$, while $\mathfrak{p}\mathcal{O}_K=\mathfrak{p}_K^e $, hence
$$\mathfrak{p}\mathcal{O}_K=\mathfrak{p}_K^{e}$$
for the unique prime (thus maximal) ideal $\mathfrak{p}_K$ of $\mathcal{O}_K$, 
and so
$$\mathcal{O}_K/\mathfrak{p}\mathcal{O}_K= \mathcal{O}_K/ \mathfrak{p}_K^{e} \cong \mathcal{O}_K/ \mathfrak{p}_K^{e}\mathcal{O}_K$$
is a finite chain ring with $q^e$ elements, if $\overline{K}=\mathbb{F}_q$. This yields an induced action of $G$ on $\mathcal{O}_K/ \mathfrak{p}_K^{e}\mathcal{O}_K$ and the above $\overline{\sigma}$ is an isomorphism of $G$-modules.
 All nonzero ideals of $\mathcal{O}_K$ are of the kind $\mathfrak{p}^l\mathcal{O}_K$ for some integer $l\geq 0$.
 Moreover, we have
 $$ \mathcal{}O_K=\mathcal{O}_F[\mu_{q-1},\varpi_K]$$

\begin{example}
(i) Let $K/F$ be the unique unramified extension of degree $n$, ${\rm Gal}(K/F)=\langle\sigma \rangle$. Then $\mathfrak{p}\mathcal{O}_K=\mathfrak{p}_K,$
$ \mathcal{O}_K/\mathfrak{p}\mathcal{O}_K= \mathcal{O}_K/ \mathfrak{p}_K=\mathbb{F}_{p^{nj}},$ and $ \overline{\sigma}\in {\rm Gal}(\mathbb{F}_{p^{nj}}/\mathbb{F}_{p^{j}}).$
\\ (ii) Let $K/F$ be a totally ramified Galois extension of degree $e$, so that $\overline{K}=\overline{F}$, then
$$\mathfrak{p}\mathcal{O}_K=\mathfrak{p}_K^{e} \text{ and }
\mathcal{O}_K/\mathfrak{p}\mathcal{O}_K= \mathcal{O}_K/ \mathfrak{p}_K^{e}\cong \mathcal{O}_K/ \mathfrak{p}_K^{e}\mathcal{O}_K$$
is a finite chain ring with $p^{je}$ elements, if $\overline{K}=\mathbb{F}_{p^j}$.
This yields an induced action of $G$ on $\mathcal{O}_K/ \mathfrak{p}_K^{e}\mathcal{O}_K$ and the above is an isomorphism of $G$-modules.
\end{example}

\subsection{Natural orders}  \label{subsec:naturalI}

Suppose
$f(t)=\sum_{i=0}^{m}d_it^i\in \mathcal{O}_K[t;\sigma,\delta]$ is monic.
 The Petit  algebra
$$S_f=K[t;\sigma,\delta]/K[t;\sigma,\delta] f$$
 has a nonassociative subalgebra given by 
$$\Lambda=\mathcal{O}_K[t;\sigma,\delta]/\mathcal{O}_K[t;\sigma,\delta] f$$
which is an order in $S_f$  (i.e., a full lattice with a nonassociative multiplicative ring structure). The order $\Lambda$ is  called the \emph{natural order} of $S_f$.
As left $\mathcal{O}_K$-module, we have  $\Lambda=\mathcal{O}_K\oplus \mathcal{O}_K t\oplus \dots\oplus \mathcal{O}_K t^{m-1}$.
If $f$ is irreducible in $K[t;\sigma,\delta]$, then $S_f$ is a division algebra and $\Lambda$ does not have any zero divisors.

The order $\Lambda$ is usually not a maximal order,
 but  it is uniquely determined whenever $Rf$ is not a two-sided ideal, since  in that case $S_f$ is a proper nonassociative algebra: $K$ is the left and middle nucleus of $S_f$ and thus uniquely determines $\mathcal{O}_{K}$ and in turn $\Lambda$.

 The center of $\Lambda$ contains  $S_0=\mathcal{O}_{F}\cap {\rm Const}(\delta)$.

 When $\delta=0$  then the  center of $\Lambda$ contains $\mathcal{O}_{F}$ and so  $\mathfrak{p}$ automatically lies in the center of $\Lambda$.
When $\delta\not=0$ we will assume from now on that the maximal ideal $\mathfrak{p}$ of $\mathcal{O}_{F}$  lies in $S_0$.
  This ensures that $\mathfrak{p}\Lambda$ is a two-sided ideal of $\Lambda$.

 For any $g(t)=\sum_{i=0}^{m-1}a_i t^i\in \mathcal{O}_K[t;\sigma,\delta]$ define $\overline{a_i}=
 a_i+\mathfrak{p}\mathcal{O}_K$ and
   $$\overline{g}(t)=\sum_{i=0}^{m-1}\overline{a_i}t^i\in (\mathcal{O}_K/\mathfrak{p}\mathcal{O}_K)[t;\overline{\sigma},\overline{\delta}].$$
  Let $\overline{f}(t)=\sum_{i=0}^{m}\overline{d_i}t^i\in (\mathcal{O}_K/\mathfrak{p}\mathcal{O}_K)[t;\overline{\sigma},\overline{\delta}]$.
  Define
  $$\Psi:\Lambda\rightarrow   S_{\overline{f}},\quad g\mapsto \overline{g}.$$
 Then $\Psi$ is a surjective homomorphism of nonassociative unital rings and ${\rm ker}(\Psi)=\mathfrak{p}\Lambda$.
For $F'=\mathbb{F}_{p^j}\cap {\rm Const}(\overline{\delta})$,
$$\Psi:\Lambda/\mathfrak{p}\Lambda\rightarrow
S_{\overline{f}},\quad g+\mathfrak{p}\Lambda \mapsto \overline{g}$$
is an $F'$-algebra isomorphism. This is proved analogously as \cite[Lemma 5]{P18.2}.

\begin{example} (analogously as in \cite[Example 6]{P18.2}) \label{ex:nca1}
Let $K/F$ be a cyclic Galois field extension of degree $n$, ${\rm Gal}(K/F)=\langle\sigma \rangle$, and $f(t)=t^n-d\in \mathcal{O}_K[t;\sigma]$ be irreducible in $K[t;\sigma]$. Then
 $S_f=(K/F,\sigma,d)$ is a  nonassociative cyclic division algebra of degree $n$ over $F$ with the natural order
 $$\Lambda=\mathcal{O}_K[t;\sigma]/\mathcal{O}_K[t;\sigma]f=\mathcal{O}_K \oplus  \mathcal{O}_K t\oplus \dots \oplus\mathcal{O}_K  t^{n-1},$$
 and
 $$\Lambda/\mathfrak{p}\Lambda\cong ((\mathcal{O}_K/\mathfrak{p}\mathcal{O}_K)/\mathbb{F}_{p^j},
\overline{\sigma},\overline{d})=S_{\overline{f}},$$
with $\overline{f}(t)=t^n-\overline{d}\in (\mathcal{O}_K/\mathfrak{p}\mathcal{O}_K)[t;\overline{\sigma}].$

If $d\in \mathcal{O}_F$ is non-zero, both the cyclic algebra $(K/F,\sigma,d)$  and the order $\Lambda$ are associative
and the choice of  $\Lambda$ depends on the choice of the maximal subfield in $(K/F,\sigma,d)$. In that case, $S_{\overline{f}}$
is a generalized associative  cyclic algebra.
If $\mathcal{O}_K/\mathfrak{p}\mathcal{O}_K$ is a field, i.e.  if $K/F$ is unramified, then $\overline{f}(t)$ is always reducible.

If $d\in \mathcal{O}_K\setminus \mathcal{O}_F$, then $(K/F,\sigma,d)$ is a proper nonassociative algebra, and  if  $1,d,\dots,d^{n-1}$ are linearly independent (e.g. when $n$ is prime) then 
$(K/F,\sigma,d)$ is a division algebra.
Here,
$$\Lambda=\mathcal{O}_K[t;\sigma]/\mathcal{O}_K[t;\sigma]f=\mathcal{O}_K \oplus  \mathcal{O}_K t\oplus \dots \oplus\mathcal{O}_K  t^{n-1}$$
 is uniquely determined.
\end{example}

\begin{remark}\label{re:I}
Suppose that $K$ is the totally unramified extension of $F$ of degree $n$. Then $\mathcal{O}_K/\mathfrak{p}\mathcal{O}_K=
\mathbb{F}_{p^{jn}}$ and  ${\rm Gal}( \mathbb{F}_{p^{jn}}/\mathbb{F}_{p^{j}} )=\langle\overline{\sigma}\rangle$.
 Let $f(t)=\sum_{i=0}^{n}d_it^i\in \mathcal{O}_K[t;\sigma,\delta]$ be such that ${\overline{f}}(t)=t^n-\overline{ d_0}$.
 \\ (i)
If $f(t)=t^n-d_0$ with $d_0\in \mathcal{O}_F$
  then the ideal
 $(\mathcal{O}_K/\mathfrak{p}\mathcal{O}_K)[t;\overline{\sigma}]\overline{f}$ is
 two-sided and the resulting $\mathbb{F}_{p^j}$-algebra associative. In this case,
$\overline{f}(t)=t^n-\overline{d_0}$ is always reducible in $\mathbb{F}_{p^{jn}}[t;\overline{\sigma}]$.
\\ (ii) If $\overline{f}$ is irreducible, then $S_{\overline{f}}$ has no non-trivial left ideals.
\end{remark}

Since we need that  $\mathfrak{p}$ lies in $S_0$, in the following we will usually assume that $\delta=0$, all arguments can be done more general, however.

\section{How to obtain $\mathbb{Z}_p$-lattice and codes from cyclic $(f,\sigma)$-codes over $\mathcal{O}_K/\mathfrak{p}\mathcal{O}_K$ using Construction A}
\label{sec:codes}

\subsection{Principal left ideals in $S_{\overline{f}}$ and $(f,\sigma)$-codes}

We keep the assumptions and notation from Section \ref{sec:naturalI} but assume that $\delta=0$.

Let  $\mathcal{I}=\Lambda g(t)$ be a principal left ideal of $\Lambda$   generated by a monic polynomial $g(t)$ such that
$\mathfrak{p}\subset \mathcal{I}$.
Then $\mathcal{I}/\mathfrak{p}\Lambda$ is a  non-zero  principal left ideal of $\Lambda/\mathfrak{p}\Lambda$ generated by the monic polynomial $g+\mathfrak{p}\Lambda$, and $\Psi(\mathcal{I}/\mathfrak{p}\Lambda)$
is a principal left ideal  of $S_{\overline{f}}$
generated by the monic polynomial $\Psi(g+\mathfrak{p}\Lambda)=\overline{ g}$. Thus the image  $\Psi(\mathcal{I}/\mathfrak{p}\Lambda)$ for any such $\mathcal{I}$ corresponds to a cyclic
$(\overline{f},\overline{\sigma})$-code $\mathcal{C}$ over $\mathcal{O}_K/\mathfrak{p}\mathcal{O}_K$, cf. Section (cf. \ref{subsec:codes}). Note that when $\mathcal{O}_K/\mathfrak{p}\mathcal{O}_K $ is a field, there is even a 1-1 correspondence between the left ideals of the algebra
$S_{\overline{f}}$ that are generated by a monic polynomial, and the $(\overline{f},\overline{\sigma})$-codes $\mathcal{C}$ over $\mathcal{O}_K/\mathfrak{p}\mathcal{O}_K$ as described in Section \ref{subsec:codes}.

We can rephrase this saying that $\Psi(\mathcal{I}/\mathfrak{p}\Lambda)$ is an $\mathcal{O}_K/\mathfrak{p}\mathcal{O}_K $-module which is isomorphic to a submodule of $\mathcal{O}_K/\mathfrak{p}\mathcal{O}_K^m$ and corresponds to a cyclic $(\overline{f},\overline{\sigma})$-code of length $m$ and dimension $m-{\rm deg}(\overline{g})$ as in \cite[Theorem  1]{BL13}.

When $\overline{f}$ is reducible and $\mathcal{O}_K/\mathfrak{p}\mathcal{O}_K$ a field, a cyclic
$(\overline{f},\overline{\sigma})$-code $\mathcal{C}$ corresponds to a right divisor $\overline{g}$ of
$\overline{f}$ and has dimension  $n-{\rm deg}(\overline{g})$.

When $\overline{f}$ is irreducible and $\mathcal{O}_K/\mathfrak{p}\mathcal{O}_K$ is a field,
then $S_{\overline{f}}$   has no nontrivial principal left ideals which contain a non-zero  polynomial
of minimal degree with invertible leading coefficient. Hence any cyclic $(\overline{f},\overline{\sigma})$-code $\mathcal{C}$ over the field $\mathcal{O}_K/\mathfrak{p}\mathcal{O}_K$ must have length $m$ and dimension $m$,
or is zero.  So we look for irreducible $f$ where
$\overline{f}$ is reducible.

In particular, when $K/F$ is a cyclic Galois extension with nontrivial automorphism $\sigma$ and we choose $f(t)$ such that $\overline{f}(t)=t^n-\overline{ c}$ with $\overline{ c}$ non-zero, then $\Psi(\mathcal{I}/\mathfrak{p}\Lambda)$
is a principal left ideal of $((\mathcal{O}_K/\mathfrak{p}\mathcal{O}_K)/\mathbb{F}_{p^j},\overline{\sigma},\overline{c})$,
 and so $\Psi(\mathcal{I}/\mathfrak{p}\Lambda)$ corresponds to a $\overline{\sigma}$-constacyclic
code.

\subsection{Construction A} \label{ex:nca}

Let  $\mathcal{I}$ be a principal left ideal of $\Lambda$ generated by a monic polynomial such that
$\mathfrak{p}\subset \mathcal{I}$. Let $\mathcal{C}$
 be the cyclic $(\overline{f},\overline{\sigma})$-code over $\mathcal{O}_K/\mathfrak{p}\mathcal{O}_K$ that corresponds to the principal left ideal $\Psi(\mathcal{I}/\mathfrak{p}\Lambda)$. Define
$$\rho: \Lambda\rightarrow  \Lambda/\mathfrak{p}\Lambda \rightarrow  \Psi(\Lambda/\mathfrak{p}\Lambda)$$
as the
canonical projection $\Lambda\rightarrow  \Lambda/\mathfrak{p}\Lambda$ composed with $\Psi$.
Then
$$L=\rho^{-1}(\mathcal{C})$$ 
is an $\mathcal{O}_F$-lattice. Its embedding into $\mathbb{Q}_p^M$ is canonically determined by considering the algebra $S_f$.
Since $\mathcal{O}_K$ is a free $\mathcal{O}_F$-module of finite rank $[\mathcal{O}_K:\mathcal{O}_F]$,
$L$ is an $\mathcal{O}_F$-module of finite rank $M=m[\mathcal{O}_K:\mathcal{O}_F]$.

When we work with $p$-adic fields, then $\mathcal{O}_F$ is a free $\mathbb{Z}_p$-module of finite rank $[\mathcal{O}_F:\mathbb{Z}_p]$ and
$L=\rho^{-1}(\mathcal{C})$ 
is a $\mathbb{Z}_p$-module of finite rank $N=m[\mathcal{O}_K:\mathbb{Z}_p]$.
The embedding of this $\mathbb{Z}_p$-lattice into $\mathbb{Q}_p^N$ is again canonically determined by
 considering $S_f$.

When the fields are Laurent series fields over finite fields and $F=\F_q(\!(t)\!)$, then $L=\rho^{-1}(\mathcal{C})$
is a $\F_q[\![t]\!]$-lattice of finite rank $N=m[\mathcal{O}_K: \F_q[\![t]\!]  ]$.

 Now all works exactly
 as explained in \cite[Section 3.3]{DO}.
The construction of the lattice $L$ can  be seen as a non-commutative and nonassociative variation of the classical Construction A in \cite{CS}.

This way we can construct a $\mathbb{Z}_p$-lattice $L$ in $\mathbb{Q}_p^N$ from the linear cyclic $(\overline{f},\overline{\sigma})$-code $\mathcal{C}$ over the finite ring
$S=\mathcal{O}_K/\mathfrak{p}\mathcal{O}_K$, and an $\F_q[\![t]\!]$-lattice $L$  in $\F_q(\!(t)\!)^N$ from the linear cyclic $(\overline{f},\overline{\sigma})$-code $\mathcal{C}$ over the finite ring
$S=\mathcal{O}_K/\mathfrak{p}\mathcal{O}_K$.

If $K/F$ is cyclic with generating automorphism $\sigma$, $f(t)=t^n-c\in \mathcal{O}_F[t]$, then $S_f$ and hence $\Lambda$ and $S_{\overline{f}}$  are associative and we are in the situation studied over number fields in \cite{DO}.

\subsection{Examples of lifting $\overline{\sigma}$-constacyclic codes}

Let $F$ be a nonarchimedean local field of residual characteristic $p\neq 2$.
Then $F^\times/(F^\times)^2 = \{1, \ep, \varpi,\ep\varpi\}$
where $\ep\in \mathcal{O}_F^\times \setminus (\mathcal{O}_F^\times)^2$ is a fixed nonsquare element of valuation zero, which is a lift to $\mathcal{O}_F^\times$ of a nonsquare element of the residue field via Hensel's Lemma \ref{le:Hensel}.
Therefore there are exactly three distinct quadratic extensions of $F$:
\begin{itemize}
\item  the unique unramified quadratic extension $K_2 = F(\sqrt{\ep})$;
\item  the ramified extension $K_\varpi = F(\sqrt{\varpi})$;
\item  a second ramified extension $K_{\ep\varpi} = F(\sqrt{\ep\varpi})$.
\end{itemize}

\begin{example}  \label{ex:1}
Let $p \not \equiv 1\, {\rm mod}\, 4$, $F=\mathbb{Q}_p$, and $K=\mathbb{Q}_p(i)$ (note that $F=K$
when $p \equiv 1\, {\rm mod}\, 4$.) Then $\mathcal{O}_F=\mathbb{Z}_p$ and $\mathcal{O}_K=\mathbb{Z}_p[i]$. 
Let $\varpi=p$ and since we assume $p\not=2$, $K$ is the unique unramified extension of $\mathbb{Q}_p$. We know that $\overline{F}=\mathbb{Z}_p/\mathfrak{p}\mathbb{Z}_p=\mathbb{F}_p$ and $\overline{K}=\mathbb{Z}_p[i]/\mathfrak{p}\mathbb{Z}_p[i]=\mathbb{F}_{p^2}$.

Define $K[t;\sigma]=\mathbb{Q}_p(i)[t;\sigma]$ with $\sigma$ the complex conjugation.
Let $f(t)=t^2-bt-c\in\mathbb{Z}_p[i][t,\sigma]$ be irreducible in $\mathbb{Q}_p(i)[t;\sigma]$.
We know that $f$ is irreducible if and only if  $\sigma(z)z-bz\not=c$ for all $z\in \mathbb{Q}_p(i)$ \cite[(17)]{P66}.
 Assume that $b=0$. Then $f(t)=t^2-c$ is irreducible if $c\not \in \mathbb{Z}_p$; in this case $A=(\mathbb{Q}_p(i)/\mathbb{Q}_p,\sigma,c)$ is a proper nonassociative quaternion division algebra over $\mathbb{Q}_p$. There exists a full parametrization of these algebras \cite{NP}.
If  $c\in \mathbb{Z}_p$ then $A=(\mathbb{Q}_p(i)/\mathbb{Q}_p,\sigma,c)$ is an associative quaternion algebra and a division algebra if and only if $c\in p \mathbb{Q}_p^{\times 2}$.
We have $\overline{f}(t)=t^2-\overline{c}\in (\mathbb{Z}_p[i]/\mathfrak{p}\mathbb{Z}_p[i])[t;\overline{\sigma}]=\mathbb{F}_{p^2}[t;\overline{\sigma}]$,
  with $\overline{f}(t)=t^2$ whenever $f(t)=t^2-c \in\mathbb{Z}_p[t,\sigma]$, i.e. whenever $A$ is an associative quaternion division algebra.

Choose $\Lambda=\mathbb{Z}_p[i]\oplus \mathbb{Z}_p[i]t$ and observe that $\Lambda$ is only unique when the algebra $A=(\mathbb{Q}_p(i)/\mathbb{Q}_p,\sigma,c)$  
is a proper nonassociative quaternion algebra. By assumption,
 $$\Lambda/\mathfrak{p}\Lambda\cong (\mathbb{F}_{p^2}/\mathbb{F}_{p}, \overline{\sigma},\overline{c})$$
 is a nonassociative quaternion algebra over $\mathbb{F}_{p}$ and
$\Lambda/\mathfrak{p}\Lambda=\mathbb{F}_{p^2}\oplus \mathbb{F}_{p^2} t$
as $\mathbb{F}_{p^2}$-vector space. This implies that for all $\overline{c}\not\in \mathbb{F}_{p}$, the nonassociative quaternion algebra $\Lambda/\mathfrak{p}\Lambda\cong (\mathbb{F}_{p^2}/\mathbb{F}_{p}, \overline{\sigma},\overline{c})$
 is a proper nonassociative  algebra over $\mathbb{F}_{p}$ and hence a division algebra. Thus in this case $\overline{f}$ is irreducible.
  Hence $\Psi(\mathcal{I}/\mathfrak{p}\Lambda)$ is either trivial or all of
$(\mathbb{F}_{p^{2}}/\mathbb{F}_{p},\overline{\sigma},\overline{c})$.

For every $c\in\mathbb{Z}_p[i]\setminus \mathbb{Z}_p$ such that  $c\in\mathfrak{p}\mathbb{Z}_p[i]=\mathfrak{p}_K$,  $\overline{f}(t)=t^2$ is reducible. In this case, $S_{\overline{f}}$ is a commutative associative algebra and
 there do not exist any
$\overline{\sigma}$-constacyclic codes since here $\overline{c}=0$. Thus this algebra
cannot be used for lattice encoding of $\overline{\sigma}$-constacyclic codes.

For every $c\in\mathbb{Z}_p[i]\setminus \mathbb{Z}_p$ such that $c\not\in\mathfrak{p}\mathbb{Z}_p[i]=\mathfrak{p}_K$, $\overline{f}(t)=t^2-\overline{c}$ is irreducible. In this case, $\Lambda/\mathfrak{p}\Lambda\cong (\mathbb{F}_{p^2}/\mathbb{F}_{p}, \overline{\sigma},\overline{c})$
 is a proper nonassociative  algebra over $\mathbb{F}_{p}$ and hence a division algebra.

For every $c\in \mathbb{Z}_p$,  $c\not\in\mathfrak{p}$ 
and
 $\overline{f}(t)=t^2-\overline{c}\in\mathcal{O}_K[t,\sigma]$ is reducible and $\overline{c}\not=0$; here  $\Lambda/\mathfrak{p}\Lambda\cong M_2(\mathbb{F}_{p})$.
For any choice of $c\in \mathbb{Z}_p$ such that $c\in\mathfrak{p}$,  $\overline{f}(t)=t^2$ is reducible as well.
\end{example}

 For a nonassociative cyclic algebra $A=(K/F,\sigma,c)$ of prime degree $n$, $A$ is a division algebra if and only if
$c\in K\setminus F$. Example \ref{ex:1} 
demonstrates that this poses a problem when $K/F$ is unramifield: 
we cannot find irreducible $f(t)=t^n-c$ such that $\overline{f}(t)=t^n-\overline{c}$  is reducible and $0\not=\overline{ c}$,
 since  $\overline{f}(t)$ is either irreducible, or $\overline{ c}=0$.  This is not the case when $n$ is not prime and when $K/F$ has ramification index $e>1$ and $f>1$.

 We conclude that proper nonassociative cyclic division algebras of prime degree $n$ over $F$ employing $K_n[t;\sigma]$ can only be used to lift $\overline{\sigma}$-constacyclic codes  of length $n$ and dimension $n$, as the arguments in the above example hold in general for $n$ prime and not just for $n=2$, hence are of less interest.

The subsequent example shows that division algebras employing $K[t;\sigma]$, where $K/F $ is a totally ramified extension,  can be used to lift $\overline{\sigma}$-constacyclic codes if and only if $c\not\in \mathfrak{p}_K$, as the arguments hold in general and not just for $n=2$:

\begin{example} \label{ex:2}

Let $K/F$ be a quadratic Galois field extension that is (totally) ramified,  then $\overline{K}=\overline{F}$ and $\overline{\sigma}=id$, e.g. choose $K=F(\sqrt{\varpi})$. Then
$\mathcal{O}_K=\mathcal{O}_F[\varpi]$,
$$\mathfrak{p}\mathcal{O}_K=\mathfrak{p}_K^{2} \text{ and }
\mathcal{O}_K/\mathfrak{p}\mathcal{O}_K= \mathcal{O}_K/ \mathfrak{p}_K^{2}\cong \mathcal{O}_K/ \mathfrak{p}_K^{2}\mathcal{O}_K.$$
is a finite chain ring with $q^2$ elements, if $\overline{K}=\mathbb{F}_q$.

Let $F=\mathbb{Q}_p$, $f(t)=t^2-c\in \mathbb{Z}_p[\varpi][t;\sigma]$, $c\in \mathbb{Z}_p[\varpi]\setminus \mathbb{Z}_p$, so that
$A=(\mathbb{Q}_p(\sqrt{\varpi})/\mathbb{Q}_p, \sigma,c)$ is a proper nonassociative quaternion division algebra over $\mathbb{Q}_p$.

For the natural order $\Lambda=\mathbb{Z}_p[\varpi]\oplus \mathbb{Z}_p[\varpi] t$ we obtain the  algebra
$$\Lambda/\mathfrak{p}\Lambda\cong ((\mathbb{Z}_p[\varpi]/ \mathfrak{p}_K^{2}\mathbb{Z}_p[\varpi] ) \mathbb{Z}_p[\varpi]/\mathbb{F}_{p},\overline{\sigma},\overline{c})$$
over  $\mathbb{F}_p$. Note that the  ring $S=\mathbb{Z}_p[\varpi]/ \mathfrak{p}_K^{2} $  has zero divisors and thus this algebra is not a division algebra independent of the choice of $c$. When $\overline{c}\not=0$ and $\overline{f}(t)=t^2-\bar c$ is reducible in $\mathbb{Z}_p[\varpi]/ \mathfrak{p}_K^{2}[t:\overline{\sigma}]$, this algebra can be used to lift codes. 

We have $\overline{c}\not=0$ if and only if $c\not\in \mathfrak{p}_K$.  If $\overline{ c}=0$ then $\overline{f}(t)=t^2$  and
$\Lambda/\mathfrak{p}\Lambda$ is a commutative associative algebra that cannot be used in our construction.
\end{example}

\section{How to obtain lattice  codes from cyclic $(\overline{f},\overline{\sigma})$-codes over $\mathcal{O}_K/\mathfrak{p}\mathcal{O}_K$ using Construction A and generalized nonassociative cyclic algebras} \label{sec:naturalII}

\subsection{Natural orders in generalized nonassociative cyclic algebras and their quotients by a prime ideal}\label{subsec:iteratedI}

Let $K/F$ be a cyclic Galois field extension of local nonarchimedean fields of degree $n$ with ${\rm Gal}(K/F)=\langle \rho\rangle$, and let
 $D=(K/F, \rho, c)$ be a cyclic division algebra over $F$ of degree $n$ such that $c\in \mathcal{O}_F^\times$. (For instance, we could choose $K=K_n$, the unique unramified extension of $F$, and $c=\varpi^m$, $1<m<n$, a power of the uniformizer.) Let $\mathcal{D}=(\mathcal{O}_K/\mathcal{O}_F,\rho, c)$
such that $\mathcal{D} \otimes_{\mathcal{O}_F}F=(K/F, \rho, c)=D$. Then
$\mathcal{D}=\mathcal{O}_K \oplus   \mathcal{O}_Ke \oplus \dots  \oplus\mathcal{O}_K  e^{n-1}$
is a natural $\mathcal{O}_F$-order of $D$, cf. \ref{subsec:naturalI} or \cite{DO}.

  Let $\sigma\in {\rm Aut}(D)$ such that
$\sigma|_{F}$ has finite order $m$. Assume that $F_0=F\cap {\rm Fix}(\sigma)$ is a local nonarchimedean field and that $F/F_0$ is an unramified Galois extension. 
Assume that $\sigma(\mathcal{D})\subset \mathcal{D}$.

Let $f(t)=t^m-d\in \mathcal{D}[t;\sigma]$. The algebra
$S_f=D[t;\sigma]/D[t;\sigma]f$
 is an example of a \emph{(generalized) nonassociative cyclic algebra of degree $m$} over $F_0$.  We denote this algebra by $(D,\sigma, d)$ and call $1,e,\dots,e^n,t,et\dots ,e^{n-1}t,$ $\dots e^nt^{m-1}$  its \emph{canonical basis} as a left $K$-vector space.
In this section we will work with generalized nonassociative cyclic algebras.

Suppose $f(t)=t^m-d\in \mathcal{D}[t;\sigma,\delta]$ is irreducible in
$D[t;\sigma,\delta]$. Consider the division algebra
$$S_f=D[t;\sigma]/D[t;\sigma] f$$
 over $F_0$. Then the $\mathcal{O}_{F_0}$-algebra
$$\Lambda=\mathcal{D}[t;\sigma]/\mathcal{D}[t;\sigma] f$$
 is an $\mathcal{O}_{F_0}$-order in $S_f$ which we call the \emph{natural order} and has no zero divisors.
 The order $\Lambda$  is not uniquely determined and
depends on the choice of the maximal subfield $K$ in $D$.

Let $1,e,\dots,e^{n-1}$ denote the canonical basis of $D$, then
$$\Lambda=\mathcal{O}_K \oplus  \mathcal{O}_K e\oplus \dots  \oplus \mathcal{O}_K e^{n-1}
\oplus \mathcal{O}_K t\oplus \mathcal{O}_K et\oplus\dots
\oplus \mathcal{O}_K e^{n-1}t \oplus \dots\oplus\mathcal{O}_K  e^{n-1}t^{m-1}$$
as a left $\mathcal{O}_K$-module.

Let $\mathfrak{p}$ denote the unique maximal ideal of $\mathcal{O}_{F_0}$. Since $F/F_0$ is unramified by our assumption, we know that
 $\mathfrak{p}\mathcal{O}_F=\mathfrak{p}_{F}$ is maximal as well.

Now $ \mathfrak{p}\mathcal{D}$ is a two-sided ideal of $\mathcal{D}$
and $\mathfrak{p}\Lambda$ is a two-sided ideal of $\Lambda$.
Let $\pi:\mathcal{D}\rightarrow  \mathcal{D}/\mathfrak{p}\mathcal{D}$ be
the canonical projection.
We have $\sigma( \mathfrak{p}\mathcal{D})\subset \mathfrak{p}\mathcal{D}$ since
 $\mathfrak{p}\subset{\rm Fix}(\sigma)$ and $\sigma(\mathcal{D})\subset \mathcal{D}$ by  assumption.
 Thus $\sigma$ induces a ring homomorphism
$$\overline{\sigma}:\mathcal{D}/\mathfrak{p}\mathcal{D} \rightarrow  \mathcal{D}/\mathfrak{p}\mathcal{D}$$
with
${\rm Fix}(\overline{\sigma})={\rm Fix}(\sigma)/\mathfrak{p}{\rm Fix}(\sigma)$
and $\pi\circ \sigma=\overline{\sigma} \circ \pi$.

Let 
$\overline{F}=\mathcal{O}_F/\mathfrak{p}\mathcal{O}_F=\mathbb{F}_{p^j}$.
 For any $g(t)=\sum_{i=0}^{m-1}a_i t^i\in \mathcal{D}[t;\sigma]$ define
   $$\overline{g}(t)=\sum_{i=0}^{m-1}\overline{a_i}t^i\in
 (\mathcal{D}/\mathfrak{p}\mathcal{D})[t;\overline{\sigma}],$$
  with $\overline{a_i}=
 a_i+\mathfrak{p},$ 
 then
  $$S_{\overline{f}}=(\mathcal{D}/\mathfrak{p}\mathcal{D})[t;\overline{\sigma}]/(\mathcal{D}/\mathfrak{p}\mathcal{D})
[t;\overline{\sigma}]\overline{f}$$
is a nonassociative ring and
an algebra over the subfield $\overline{\mathcal{O}_{F_0}}$ of $\mathbb{F}_{p^j}$. Then
$$\Psi:\Lambda\rightarrow   (\mathcal{D}/\mathfrak{p}\mathcal{D})[t;\overline{\sigma}]/
(\mathcal{D}/\mathfrak{p}\mathcal{D})[t;\overline{\sigma}]\overline{f},\quad g\mapsto \overline{g}$$
is a  surjective homomorphism of additive groups with kernel $\mathfrak{p}\Lambda$, and
induces an $\overline{\mathcal{O}_{F_0}}$-algebra isomorphism
$$\Psi:\Lambda/\mathfrak{p}\Lambda\cong
(\mathcal{D}/\mathfrak{p}\mathcal{D})[t;\overline{\sigma}]/(\mathcal{D}/\mathfrak{p}\mathcal{D})
[t;\overline{\sigma}]\overline{f},\quad g+\mathfrak{p}\Lambda \mapsto \overline{g}$$
 (The proof is analogous to the one of \cite[Lemma 12]{P18.2}.)

\subsection{Lattice encoding of cyclic $(f,\sigma)$-codes}\label{sub:ex:nca}

Let $\mathcal{I}$ be a  principal left ideal of $\Lambda$ generated by a monic polynomial $g(t)$, such that
$\mathfrak{p}\subset \mathcal{I}\cap \mathcal{O}_{F_0}$.
Then $\mathcal{I}/\mathfrak{p}\Lambda$ is a non-zero principal left ideal of $\Lambda/\mathfrak{p}\Lambda$ and
$\Psi(\mathcal{I}/\mathfrak{p}\Lambda)$
is a principal left ideal of $S_{\overline{f}}$ generated by the monic polynomial $\Psi(g+\mathfrak{p}\Lambda)=\overline{g}$.
Again, $\Psi(\mathcal{I}/\mathfrak{p}\Lambda)$ corresponds to a cyclic
$(\overline{f},\overline{\sigma})$-code $\mathcal{C}$ over $\mathcal{O}_K/\mathfrak{p}\mathcal{O}_K$.

If  $f(t)$ satisfies
$\overline{f}(t)=t^n-\overline{ c}\in \mathcal{D}/\mathfrak{p}\mathcal{D}$ with $\overline{c}$ non-zero, then $\Psi(\mathcal{I}/\mathfrak{p}\Lambda)$ is a
$\overline{\sigma}$-constacyclic code over $\mathcal{O}_K/\mathfrak{p}\mathcal{O}_K$.

If $\overline{f}$ is irreducible and $\mathcal{D}/\mathfrak{p}\mathcal{D}$ is a division algebra, then the algebra
$S_{\overline{f}}$ is simple. Then any non-zero code
$\mathcal{C}$ must have length $m$ and dimension $m$ (and correspond to the whole algebra),
whereas whenever $\overline{f}$ is reducible, $\mathcal{C}$ respectively $\Psi(\mathcal{I}/\mathfrak{p}\Lambda)$
 corresponds to a monic right divisor $\overline{g}$ of
$\overline{f}$ and has dimension  $m-{\rm deg}(\overline{g})$.
Let
$$\rho: \Lambda\rightarrow  \Lambda/\mathfrak{p}\Lambda \rightarrow  \Psi(\Lambda/\mathfrak{p}\Lambda)$$ be the
canonical projection  composed with $\Psi$.
We know that $\mathcal{O}_K$ is a free $\mathcal{O}_{F_0}$-module of rank $[\mathcal{O}_K:\mathcal{O}_{F_0}]$, and
a free $\mathbb{Z}_p$-module of rank $[\mathcal{O}_K:\mathbb{Z}_p]$ when we work with $p$-adic fields.
 Again all works analogously as in \cite[Section 3.3]{DO}: When the fields are $p$-adic then
$$L=\rho^{-1}(\mathcal{C})$$
is a $\mathbb{Z}_p$-lattice which is a $\mathbb{Z}_p$-module of dimension $N=n^2 m[\mathcal{O}_F:\mathbb{Z}_p]$.
The embedding of this $\mathbb{Z}_p$-lattice into $\mathbb{Q}_p^N$ is canonically determined by
 considering $S_f$.

When the fields are Laurent series over finite fields with $F_0=\F_q(\!(t)\!)$, then
$L=\rho^{-1}(\mathcal{C})$
is a $\F_q[\![t]\!]$-lattice which is a $\F_q[\![t]\!]$-module of dimension $N=n^2 m[\mathcal{O}_F:\F_q[\![t]\!]]$.

This is another generalization of  Construction A.

\section{Lifting linear matrix codes to maximum rank distance  codes}\label{sec:MRD}

\subsection{Rank metric codes}

Let $K$ be a field. In the context of rank metric codes, a \emph{code} is a set  of matrices $\mathcal{C}\subset M_{n\times m}(K)$.
Let $L\subset K$ be a subfield. A code $\mathcal{C}\subset M_{n\times m}(K)$ is \emph{$L$-linear}, if $\mathcal{C}$ is an $L$-vector space.
A \emph{rank metric code} is a code $\mathcal{C}\subset M_{n, m}(K)$
 equipped with the rank distance function $d(A,C)=\text{rank}(A-C).$ The \textit{minimum distance} of a rank metric code $\mathcal{C}$ is defined to be
 $$d_{\mathcal{C}}=\text{min}\lbrace d(A,C)\mid  A,C\in \mathcal{C}, A\neq C\rbrace.$$
 An $L$-linear rank metric code $\mathcal{C}$ satisfies the Singleton-like bound
  $$dim_{L}(\mathcal{C})\leq n(m-d_{\mathcal{C}}+1)[K:L],$$
where $dim_L(\mathcal{C})$ is the dimension of the $L$-vector space $\mathcal{C}$ \cite[Proposition 6]{ALR}.

An $L$-linear rank metric code attaining the Singleton-like bound is called a  \textit{maximum rank distance code} or \textit{MRD-code} (for MRD-codes over cyclic field extensions see \cite{ALR}).

More generally, we can also say that a \emph{code} is a set  of matrices $\mathcal{C}\subset M_{n\times  m}(B)$, where $B$ is a division algebra.
Let $B'\subset B$ be a subalgebra, then $\mathcal{C}$ is \emph{$B'$-linear} (or simply \emph{linear}) if $\mathcal{C}$ is a right $B'$-module.

Let $\text{colrank}(A)$ be the column rank of a matrix $A \in M_{n\times m}(B)$ (i.e., the rank of the right $B$-module generated by the columns of $A$). A matrix in $M_{n\times m}(B)$ has column rank at most $m$; any matrix which attains this bound is said to have attained \emph{full column rank}. A \emph{generalized rank metric code} is a code $\mathcal{C}\subset M_{n \times  m}(B)$ together with the distance function defined as
$$d(A,C)=\text{colrank}(A-C),$$
for all $A,C\in M_{n\times m}(B)$.

The \emph{minimum distance} of a generalized  rank metric code $\mathcal{C}\subset M_{n\times m}(B)$ is defined as
$$d_{\mathcal{C}}=\text{min}\lbrace d(A,C)\mid  A,C\in \mathcal{C}, A\neq C\rbrace,$$
cf. \cite{PT}.

Rank distance codes are used in coding theory, e.g. for error correcting codes \cite{R91},
for designing  public-key cryptosystems, e.g. \cite{GPT91, FL06}, in distributed data storage, and in coded caching.  For an excellent overview see \cite{Betal}. Traditionally, MRD codes  are studied over finite fields. Recently they were also considered over  number fields \cite{ALR}.
Their definition has  also been generalized to codes with entries in noncommutative division rings \cite{PT}.

\subsection{General case, employing $f\in K[t;\sigma]$}
Let $K/F$ be a Galois extension of local nonarchimedean fields of degree $n$ and let  $\sigma\in G={\rm Gal}(K/F)$.
Let
$f(t)=\sum_{i=0}^{m}d_it^i\in \mathcal{O}_K[t;\sigma]$
 be monic and irreducible  in $K[t;\sigma]$, such that
$S_f=K[t;\sigma]/K[t;\sigma] f$ is a division algebra, and
$\Lambda=\mathcal{O}_K[t;\sigma]/\mathcal{O}_K[t;\sigma] f=\mathcal{O}_K\oplus \mathcal{O}_K t\oplus \dots\oplus \mathcal{O}_K t^{m-1}$ be the natural order in the Petit algebra $A=S_f$. Identify  $x = x_0 + x_1t + \cdots +  x_{m-1}t^{m-1}\in \Lambda$  with  $(x_0, x_1, \ldots, x_{m-1})\in \mathcal{O}_K^m$.
Then the multiplication in $\Lambda$ can be written as
$$(x_0,\dots,x_{m-1})\cdot  (y_0,\dots,y_{m-1})=M(x_0,\dots,x_{m-1})(y_0,\dots,y_{m-1})^t$$
and the left multiplication $L_x\in {\rm End}_{\mathcal{O}_K}(\Lambda)$  in $\Lambda$  can be represented by an invertible matrix $M(x) \in M_m(\mathcal{O}_K)$.

The family of matrices $\mathcal{C}=\{M(x)\,|\, 0\not=x\in \Lambda\}$ is a linear MRD-code over the ring of integers $\mathcal{O}_K$: since $\Lambda$ is an order in the division algebra $A$,
each matrix $M(x)$  must have maximum column rank $m$, for all $x\in \Lambda$, $x\not=0$ (assume to the contrary that this is not the case for some $x\in \Lambda$, then $M(x)$ read as a matrix over $K$ would not be invertible, a contradiction to the fact that $A$ is a division algebra).

Let $\mathfrak{p}$ be the maximal ideal of $\mathcal{O}_{F}$ and let $\mathcal{I}$ be a principal left ideal of $\Lambda$   generated by a monic polynomial such that
$\mathfrak{p}\subset \mathcal{I}$. Then $\mathcal{C}_{\mathcal{I}}=\{M(x)\,|\, 0\not=x\in \mathcal{I}\}\subset \mathcal{C}$ is also a linear MRD-code over the ring of integers $\mathcal{O}_K$; all its matrices have maximum column rank (and nonzero determinant when read as matrices over $K$).

Since $\mathcal{I}$ is a principal left ideal of $\Lambda$   generated by a monic polynomial, say $g(t)$, $\Psi(\mathcal{I}/\mathfrak{p}\Lambda)$ is principal left ideal  generated by $\overline{g}(t)$.

The multiplication in
$$\overline{A}=S_{\overline{f}}=(\mathcal{O}_K/\mathfrak{p}\mathcal{O}_K)/(\mathcal{O}_K/\mathfrak{p}\mathcal{O}_K)\overline{f}$$
 can be written as
$$(x_0,\dots,x_{m-1})\cdot  (y_0,\dots,y_{m-1})=M_0(x_0,\dots,x_{m-1})(y_0,\dots,y_{m-1})^t$$
where the matrix $M_0(x)\in M_{f_0}(\mathcal{O}_K/\mathfrak{p}\mathcal{O}_K)$ represents the left multiplication $L_x\in {\rm End}_{\overline{K}}(\overline{A})$ with $f_0=[\mathcal{O}_K/\mathfrak{p}\mathcal{O}_K:\overline{F}]$. Since $\overline{A}$ need not be a division algebra, a matrix $M_0(x)$ may have nontrivial determinant. Nonetheless, the set of matrices $\{M_0(x)\,|\, 0\not=x\in \overline{A}\}$ is a linear code over the finite ring $\mathcal{O}_K/\mathfrak{p}\mathcal{O}_K$, and so is the subset  of matrices $\mathcal{C}_0=\{M_0(x)\,|\, 0\not=x\in \Psi(\mathcal{I}/\mathfrak{p}\Lambda)\}$.

So let $\mathcal{C}_0$ be the linear matrix code over $\mathcal{O}_K/\mathfrak{p}\mathcal{O}_K$ that corresponds to the principal left ideal $\Psi(\mathcal{I}/\mathfrak{p}\Lambda)$.
Use the canonical projection composed with $\Psi$, $\rho: \Lambda\rightarrow  \Lambda/\mathfrak{p}\Lambda \rightarrow  \Psi(\Lambda/\mathfrak{p}\Lambda)=
((\mathcal{O}_K/\mathfrak{p}\mathcal{O}_K)/\mathbb{F}_{p^f},
\overline{\sigma},\overline{c})$,  to lift the linear code $\mathcal{C}_0$ to
$L=\rho^{-1}(\mathcal{C}_0)$ which is a linear MRD-code $\mathcal{C}_{\mathcal{I}}$ over $\mathcal{O}_K$ and an $\mathcal{O}_F$-lattice of rank $N=n[\mathcal{O}_K:\mathcal{O}_F]$.

When working with $p$-adic fields, $L$  is a linear MRD-code over $\mathcal{O}_K$ that
is a $\mathbb{Z}_p$-lattice of rank $N=m[\mathcal{O}_K:\mathbb{Z}_p]$. When working with fields that are Laurent series over finite fields with $F=\F_q(\!(t)\!)$, then
$L=\rho^{-1}(\mathcal{C})$
is a $\F_q[\![t]\!]$-module of dimension $N=m[\mathcal{O}_K:\F_q[\![t]\!]]$.

We only considered the case $\delta=0$ here but everything works more generally, too. We now look at the special case of lifting matrix codes obtained from nonassociative cyclic algebras:

\subsection{MRD-codes from nonassociative cyclic algebras}

 Let $A=(K/F,\sigma,c)$, $c\in\mathcal{O}_K$ non-zero, be a nonassociative cyclic division algebra  over $F$
 of degree $m$ with $c\in \mathcal{O}_K$ and let $\Lambda=\mathcal{O}_K[t;\sigma]/\mathcal{O}_K[t;\sigma]f$ (we allow the algebra to be associative).
 The left multiplication $L_x\in {\rm End}_{\mathcal{O}_K}(\Lambda)$  in $\Lambda$  is represented by the  matrix
\begin{equation} \label{equ:matrix_rep_cda}
M(x)= \left[ \begin{array}{ccccc}
x_0 & c \sigma(x_{m-1})& c \sigma^2(x_{m-2}) & \cdots & c \sigma^{m-1}(x_1) \\
x_1 & \sigma(x_0) & c\sigma^2(x_{m-1}) & \cdots & c \sigma^{m-1}(x_{2}) \\
x_2 & \sigma(x_1) & \sigma^2(a_0) & \cdots & c \sigma^{m-1}(x_3)\\
\vdots & \vdots & \vdots & \ddots & \vdots \\
x_{m-1} & \sigma(x_{m-2}) & \sigma^2(x_{m-3}) & \cdots & \sigma^{m-1}(x_0) \end{array} \right]  \in M_m(\mathcal{O}_K)
\end{equation}
which has column rank $m$, since $A$ is a division algebra.
The family of matrices $\mathcal{C}=\{M(x)\,|\, 0\not=x\in \Lambda\}$ is thus a linear MRD-code over the ring of integers $\mathcal{O}_K$. 

Let $\mathfrak{p}$ be the maximal ideal of $\mathcal{O}_F$  and let $\mathcal{I}$ be a principal left ideal of $\Lambda$  such that $\mathfrak{p}\subset \mathcal{I}$ that is generated by a monic polynomial. Then $\mathcal{C}_{\mathcal{I}}=\{ M(x) \,|\, 0\not=x\in \mathcal{I}  \}$ is a linear MRD-code over the ring of integers $\mathcal{O}_K$.

For
$$\overline{A}=( (\mathcal{O}_K/\mathfrak{p}\mathcal{O}_K) / \overline{F}, \overline{\sigma}, \overline{c}),$$
 the left multiplication $L_x \in {\rm End}_{\overline{K}}(\overline{A})$ is represented by the matrix
\begin{equation} \label{equ:matrix_rep_cda}
M_0(x)= \left[ \begin{array}{ccccc}
x_0 & \overline{c} \sigma(x_{f-1})& \overline{c} \sigma^2(x_{f-2}) & \cdots & \overline{c} \sigma^{f-1}(x_1) \\
x_1 & \sigma(x_0) & \overline{c}\sigma^2(x_{f-1}) & \cdots & \overline{c} \sigma^{m-1}(x_{2}) \\
x_2 & \sigma(x_1) & \sigma^2(a_0) & \cdots & \overline{c} \sigma^{f-1}(x_3)\\
\vdots & \vdots & \vdots & \ddots & \vdots \\
x_{m-1} & \sigma(x_{f-2}) & \sigma^2(x_{f-3}) & \cdots & \sigma^{f-1}(x_0) \end{array} \right]  \in M_f(\mathcal{O}_K/\mathfrak{p}\mathcal{O}_K)
\end{equation}
with $f=[\mathcal{O}_K/\mathfrak{p}\mathcal{O}_K:\overline{F}]$  and the set  of matrices $\mathcal{C}_0=\{M_0(x)\,|\, 0\not=x\in \Psi(\mathcal{I}/\mathfrak{p}\Lambda)\}$ is a linear code over $\mathcal{O}_K/\mathfrak{p}\mathcal{O}_K$.

 So let $\mathcal{C}_0$ be the linear matrix code over $\mathcal{O}_K/\mathfrak{p}\mathcal{O}_K$ that corresponds to the principal left ideal $\Psi(\mathcal{I}/\mathfrak{p}\Lambda)$. Then
 $L=\rho^{-1}(\mathcal{C}_0)$ is a linear MRD-code over $\mathcal{O}_K$ and an $\mathcal{O}_F$-lattice of rank $N=n[\mathcal{O}_K:\mathcal{O}_F]$.

When working with $p$-adic fields, $L$  is a linear MRD-code $\mathcal{C}_{\mathcal{I}}$ over $\mathcal{O}_K$ that
is a $\mathbb{Z}_p$-lattice of rank $N=m[\mathcal{O}_K:\mathbb{Z}_p]$ as $\mathbb{Z}_p$-module. When working with fields that are Laurent series over finite fields with $F=\F_q(\!(t)\!)$, then
$L=\rho^{-1}(\mathcal{C})$
is an $\F_q[\![t]\!]$-lattice of rank $N=m[\mathcal{O}_K:\F_q[\![t]\!]]$ as $\F_q[\![t]\!]$-module.

\subsection{MRD-codes from nonassociative generalized cyclic algebras} \label{ex:gencyclic}
Let $F$, $L$ and $F_0=F\cap L$ be local nonarchimedean fields, and let $K$ be a cyclic field extension of both $F$ and $L$ such that $[K:F] = n$, $[K:L] = m$, and ${\rm Gal}(K/F) = \langle \rho\rangle$, ${\rm Gal}(K/L) = \langle \sigma \rangle$. Suppose that $F/F_0$ is unramified.
  We must assume that $\rho\circ \sigma=\sigma\circ \rho$.

 Let $D=(K/F, \rho, c)$  be an associative cyclic division algebra over $F$  with canonical basis
 $1,e,\dots,e^{n-1}$
 and  $c\in F_0$.
 For $x= x_0 + x_1 e+x_2  e^2+\dots + x_{n-1}e^{n-1}\in D$, extend $\sigma$ to an automorphism $\sigma\in {\rm Aut}_L(D)$
 of order $m$  via
$$\sigma(x)=\sigma(x_0) +  \sigma(x_1)e +\sigma(x_2)  e^2 +\dots +\sigma(x_{n-1}) e^{n-1}.$$
 For all $d \in D^\times$, the generalized nonassociative cyclic algebra
$S_f=D[t;\sigma]/D[t;\sigma](t^m-d)=(D,\sigma,d)$
has dimension $m^2n^2$ over $F_0$. For all $d \in F^\times$, we have
$(D, \sigma, d)= (L/F_0,\rho, c)\otimes_{F_0} (F/F_0, \sigma, d)$ and the algebra
 is associative if and only if $d\in F_0$.

 Let $c \in \mathcal{O}_{F_0}$ then $\mathcal{D}=(\mathcal{O}_K/\mathcal{O}_F, \rho, c)$ has degree $n=[K:F]$ and
 $\mathcal{D} \otimes_{\mathcal{O}_F}F=D$.
 For $x= x_0 + x_1 e+x_2 e^2 +\dots + x_{n-1}e^{n-1}\in D$, define
$$\sigma(x)=\sigma(x_0) +  \sigma(x_1)e + \sigma(x_2)e^2 +\dots + \sigma(x_{n-1})e^{n-1}.$$
Since $c \in \mathcal{O}_{F_0}$, $\sigma\in {\rm Aut}_L(D)$ has order $m$
and restricts to $\sigma\in {\rm Aut}_{\mathcal{O}_{L}}(\mathcal{D})$.

 Let $\mathfrak{p}$ be the maximal ideal in $\mathcal{O}_{F_0}$ then by our assumption that $F/F_0$ is unramified, $\mathfrak{p}\mathcal{O}_F=\mathfrak{p}_K$ is maximal, too.  Let $d\in \mathcal{O}_L^\times$ or
$d\in \mathcal{O}_F^\times$ and define   $A=(D,\sigma,d)$. Then
$$\Lambda/\mathfrak{p}\Lambda\cong
(\mathcal{D}/\mathfrak{p}\mathcal{D})[t;\overline{\sigma}]/(\mathcal{D}/\mathfrak{p}\mathcal{D})
[t;\overline{\sigma}]\overline{f}\cong (\overline{D}, \overline{\sigma}, \overline{d})$$
is an algebra over
 $\overline{F_0}=\mathcal{O}_{F_0}/\mathfrak{p}$,
with
$\overline{D}=\mathcal{D}/\mathfrak{p}\mathcal{D}$
 a generalized associative cyclic algebra over $\mathbb{F}_{p^j}={\rm Fix}(\overline{\rho})$.

 Assume the situation of Section \ref{ex:gencyclic}, i.e. $A=(D,\sigma,d)$ be a division algebra of degree $n$ and $d\in \mathcal{O}_L$ or $d\in \mathcal{O}_F$.

For $x=x_0+x_1t+x_2t^2+\cdots+x_{m-1}t^{m-1}$, $y=y_0+y_1t+y_2t^2+\cdots+y_{m-1}t^{m-1}\in \Lambda$ with $x_{i}, y_i \in \mathcal{D}$, represent $x$ as $(x_0, x_1, \ldots, x_{m-1})$, $y$ as
$(y_0, y_1, \ldots, y_{m-1})$, and write the multiplication in $\Lambda$  as
$$(x_0,\dots,x_{m-1})\cdot  (y_0,\dots,y_{m-1})=M'(x_0,\dots,x_{m-1})(y_0,\dots,y_{m-1})^t$$
with
\[M'(x) = \left[ \begin{array}{ccccc}
x_0 & d \sigma(x_{m-1})& d \sigma^{2}(x_{m-2}) & \cdots & d \sigma^{m-1}(x_1) \\
x_1 & \sigma(x_0) & d \sigma^{2}(x_{m-1}) & \cdots & d \sigma^{m-1}(x_{2}) \\
x_2 & \sigma(x_1) & \sigma^{2}(x_0) & \cdots & d \sigma^{m-1}(x_3)\\
\vdots & \vdots & \vdots & \ddots & \vdots \\
x_{m-1} & \sigma(x_{m-2}) & \sigma^{2}(x_{m-3}) & \cdots & \sigma^{m-1}(x_0) \end{array} \right] \in M_m(\mathcal{D}). \]
 Substitute every entry $d$ in $M'(x)$ with the left regular representation $\gamma(d)$ in $\mathcal{D}$, and  every
entry $x_i$ in $M'(x)$ with the left regular representation $\gamma(x_i)$ in $\mathcal{D}$. This yields a matrix
\[M(x)=\gamma(M'(x)) = \left[ \begin{array}{cccc}
\gamma(x_0) & \gamma(d) \sigma(\gamma(x_{m-1}))&  \cdots & \gamma(d) \sigma^{m-1}(\gamma(x_1)) \\
\gamma(x_1) & \sigma(\gamma(x_0)) & \cdots & \gamma(d) \sigma^{m-1}(\gamma(x_{2})) \\
                                       \vdots & \vdots  & \ddots & \vdots \\
\gamma(x_{m-1}) & \sigma(\gamma(x_{m-2})) & \cdots & \sigma^{m-1}(\gamma(x_0)) \end{array} \right]\in M_{mn}(\mathcal{O}_K)
\]
where $\sigma(\gamma(x_i))$ means we apply $\sigma$ to each entry of the matrix $\gamma(x_i)\in M_n (\mathcal{O}_K)$.
Products are the usual matrix products. This
matrix represents left multiplication $L_x\in {\rm End}_\mathcal{D}(\Lambda)$ in $\Lambda$.
 Write elements in
$\Lambda=\mathcal{O}_K \oplus  \mathcal{O}_K e\oplus  \dots\oplus\mathcal{O}_K  e^{n-1}t^{m-1}$
as row vectors of length $mn$ with entries in $\mathcal{O}_K$. Write the multiplication in  $\Lambda$ as
$$(x_0,\dots,x_{mn-1})\cdot  (y_0,\dots,y_{mn-1})=M(x_0,\dots,x_{mn-1})(y_0,\dots,y_{mn-1})^t.$$
The family of matrices $\mathcal{C}=\{ M(x)\,|\, x\in \Lambda\}$ is a linear MRD-code, as all matrices have column rank $m$. In particular,
if $d\in \mathcal{O}_F$, then
$\det(\gamma(M'(x))) \in \mathcal{O}_F$ (analogously as proved in \cite{MO13}, or cf. \cite[Remark 5]{PS15})
and if $d\in \mathcal{O}_L$, then
\[ \gamma(M'(x)) = \begin{bmatrix}
             \gamma(x_0) & d \sigma(\gamma(x_{n-1}))& d \sigma^{2}(\gamma(x_{n-2})) & \cdots & d \sigma^{m-1}(\gamma(x_1)) \\
             \gamma(x_1) & \sigma(\gamma(x_0)) & d \sigma^{2}(\gamma(x_{n-1})) & \cdots & d \sigma^{m-1}(\gamma(x_{2})) \\
             \gamma(x_2) & \sigma(\gamma(x_1)) & \sigma^{2}(\gamma(x_0)) & \cdots & d \sigma^{m-1}(\gamma(x_3))\\
             \vdots & \vdots & \vdots & \ddots & \vdots \\
             \gamma(x_{n-1}) & \sigma(\gamma(x_{n-2})) & \sigma^{2}(\gamma(x_{n-3})) & \cdots & \sigma^{m-1}(\gamma(x_0))
             \end{bmatrix}
             \]
and $\det(\gamma(M(x))) \in L\cap \mathcal{O}_K=\mathcal{O}_L$ (analogously as shown in \cite{R13}, \cite[Lemma 19]{PS15}).

 Let $\mathfrak{p}\subset \mathcal{O}_{F_0}$ be the maximal  ideal then $\mathfrak{p}\mathcal{O}_F=\mathfrak{p}_F$.
 Let $\mathcal{I}$ be a principal left ideal of $\Lambda$   generated by a monic polynomial such that
$\mathfrak{p}\subset \mathcal{I}$.  Then $\mathcal{C}_\mathcal{I}=\{M(x)\,|\, 0\not=x\in \mathcal{I}\}$ is a linear MRD-code over the ring of integers $\mathcal{O}_K$, all its matrices have nonzero determinant.

 If $\mathcal{I}$ is generated by the monic polynomial $g(t)$, then $\Psi(\mathcal{I}/\mathfrak{p}\Lambda)$ is a principal left ideal  generated by $\overline{g}(t)$.

The multiplication in
$$\overline{A}=(\overline{D}, \overline{\sigma}, \overline{d})$$
 can be written as
$$(x_0,\dots,x_{mn-1})\cdot  (y_0,\dots,y_{mn-1})=M_0(x_0,\dots,x_{mn-1})(y_0,\dots,y_{mn-1})^t$$
where the left multiplication $L_x\in {\rm End}_{\overline{F}_{0}}(\overline{A})$  is represented by the matrix
 \[M_0(x) = \begin{bmatrix}
             \gamma(x_0) & \overline{d} \sigma(\gamma(x_{n-1}))& \overline{d} \sigma^{2}(\gamma(x_{n-2})) & \cdots & \overline{d} \sigma^{m-1}(\gamma(x_1)) \\
             \gamma(x_1) & \sigma(\gamma(x_0)) & \overline{d} \sigma^{2}(\gamma(x_{n-1})) & \cdots & \overline{d} \sigma^{m-1}(\gamma(x_{2})) \\
             \gamma(x_2) & \sigma(\gamma(x_1)) & \sigma^{2}(\gamma(x_0)) & \cdots & \overline{d} \sigma^{m'-1}(\gamma(x_3))\\
             \vdots & \vdots & \vdots & \ddots & \vdots \\
             \gamma(x_{n-1}) & \sigma(\gamma(x_{n-2})) & \sigma^{2}(\gamma(x_{n-3})) & \cdots & \sigma^{m-1}(\gamma(x_0))\in
             \end{bmatrix}
             \]
 Then $\mathcal{C}_0=\{M_0(x)\,|\, 0\not=x\in \Psi(\mathcal{I}/\mathfrak{p}\Lambda)\}$ is a linear code.
 Again, the pre-image of $\mathcal{C}_0$  under
$$\rho: \Lambda\rightarrow  \Lambda/\mathfrak{p}\Lambda \rightarrow  \Psi(\Lambda/\mathfrak{p}\Lambda)
=(\overline{D}, \overline{\sigma}, \overline{d}),$$
 yields a linear MRD-code $L=\rho^{-1}(\mathcal{C}_0)$ over $\mathcal{O}_K$ that is an $\mathcal{O}_{F_0}$-lattice.

When we work with $p$-adic fields, $L=\rho^{-1}(\mathcal{C})$
is a linear MRD-code over $\mathcal{O}_K$ which
is a $\mathbb{Z}_p$-lattice whose  embedding into $\mathbb{Q}_p^N$, $N=mn^2[\mathcal{O}_{F_0}:\mathbb{Z}_p]$, is canonically determined by
the multiplication in $A$.

 When we work with Laurent series over finite fields and $F_0=\F_q(\!(t)\!)$, then
$L=\rho^{-1}(\mathcal{C})$
is a $\F_q[\![t]\!]$-lattice of dimension $N=mn^2[\mathcal{O}_{F_0}:\F_q[\![t]\!]]$.

This construction also works when we choose any monic $f(t)=\sum_{i=0}^{m}d_it^i\in \mathcal{D}[t;\sigma]$ that is irreducible
in $K[t;\sigma]$. In particular it also works for associative algebras, e.g. generalized cyclic algebras that are associative.

\begin{remark}
The explanations in \cite[Section 5.2, 5.3]{DO} hold analogously for our generalizations of Construction A using rings of integers of local nonarchimedean fields, and show the potential of the construction for coset coding used in space-time block coding, in particular for wiretap
 space-time block coding, but also for linear codes over finite rings.
\end{remark}

\section{Outlook} \label{sec:outlook}

\subsection{A $p$-adic lattice approach at learning with errors}   The lack of diversity among post-quantum assumptions is widely recognized as a big problem. Working with $p$-adic lattices may be a step towards  addressing this. While the author is no expert in this area, she hopes that recent advantages in post-quantum cryptography, for example cf. \cite{Deng, Detal, Z, ZDW, Z2024} will open the door to this exciting new development.

It would be interesting to see how the  $p$-adic lattices constructed in this paper can be employed in lattice based cryptography,
 for instance in the context of recent novel variants of learning with errors (LWE).  To the author's knowledge, so far only  rings of integers of number fields or modules over these rings have been used in these constructions.
  Indeed, lattices arising from proper nonassociative cyclic algebras over number fields have just been used for the first time in \cite{MC}  (CLWE, NCLWE).

The approach looks promising as all of the mathematical results obtained in \cite{MC} on the multiplicative
ideal theory of the two-sided ideals in a nonassociative natural order $\Lambda$ analogously hold  when $K/F$ is a Galois extension of local nonarchimedean fields (and not a number field extension).   Analogues of
Minkowski's first and second theorem for Euclidean lattices for  $\mathcal{O}_F$-lattices have recently been proved for $F$ a local nonarchimedean field \cite{Deng}.

 Let us sketch our approach for the setting in this paper: Let $\Lambda$ be the natural order in a proper nonassociative cyclic division algebra $A=(K/F,\sigma, a)$ of degree $n$, $a\in \mathcal{O}_K^\times \setminus \mathcal{O}_F^\times$. Then \cite[Proposition 5, 6, 7]{MC} and \cite[Lemma 3, Theorem 5, Corollary 1]{MC} hold analogously for local nonarchimedean fields.

\begin{theorem} (This is \cite[Theorem 6]{MC} for local nonarchimedean fields)
The multiplication of $\Lambda$-ideals $\mathcal{I}$ where $\mathcal{I}\cap  \mathcal{O}_F$ is unramified in $ \mathcal{O}_K$ yields ideals and is commutative and associative.
\end{theorem}

The proof is verbatim as the proof of \cite[Theorem 6]{MC}. This result might be eventually used to  give  some variation of an unrestricted search-to-decision reduction for NCLWE samples, as over number fields.

 The same definitions of unramifield ideals in $\Lambda$, inverse ideals and  dual lattices when working with number fields and their rings of integers that are given in \cite{MC} can also be used to define inverse ideals and  dual lattices in our setting. Then \cite[Proposition 8, 9]{MC} hold analogously for $p$-adic fields, now using the trace of the field extension $K/\mathbb{Q}_p$, respectively of the extension $K/\F_q(\!(t)\!)$:

\begin{proposition} (This is \cite[Proposition 9]{MC}  for local nonarchimedean fields)
Let $\mathcal{I}\subset  \lambda$ be a two-sided integral unramified ideal.
 \\ (i) Let $F$ be a $p$-adic field. Define $\mathcal{J}_\mathcal{I}=\{x\in A\,|\, tr(xy)\in \mathbb{Z}_p \text{ for all } x\in \mathcal{I}\}$. Then the dual ideal of $\mathcal{I}$ is given by $\mathcal{J}_\mathcal{I} $.
\\ (ii)
When we work with Laurent series over finite fields and $F_0=\F_q(\!(t)\!)$, i.e. with $\F_q[\![t]\!]$-lattices, define $\mathcal{J}_\mathcal{I}=\{x\in A\,|\, tr(xy)\in \F_q[\![t]\!] \text{ for all } x\in \mathcal{I}\}$. Then the dual ideal of $\mathcal{I}$ is given by $\mathcal{J}_\mathcal{I} $.
\end{proposition}

The proof is verbatim as the proof of \cite[Proposition 9]{MC}.

Since we do not have a Frobenius norm in our setting, the next step in trying to imitate the setup presented in \cite{MC} would be to use the maximum norm defined in \cite[Definition 3.9]{Deng}.

For local nonarchimedean fields $F$, we note that the dual lattice of an $\mathcal{O}_F$-lattice, in particular of a $p$-adic lattice, is defined in \cite{Deng} as well, and analogues of
Minkowski's first and second theorem for Euclidean lattices for  $\mathcal{O}_F$-lattices  are given. The next step would be to compare these two definitions of a dual lattice.

\subsection{Orders other than natural orders} The theory of orders in associative central simple algebras is well developed.  The natural orders $\Lambda$ are generally not maximal.  An obvious next step is to employ orders in the algebras $S_f$ other than the natural orders $\Lambda$ defined here.

First attempts to classify classes of maximal orders in a proper  nonassociative quaternion algebras have been made in \cite{Ka, LW}. The results there can be generalized to proper nonassociative cyclic algebras of any degree $n$. It is clear that the flavour of the theory changes substantially for proper nonassociative quaternion algebras already. This, however, should make the resulting lattice codes more interesting.


\providecommand{\bysame}{\leavevmode\hbox to3em{\hrulefill}\thinspace}
\providecommand{\MR}{\relax\ifhmode\unskip\space\fi MR }
\providecommand{\MRhref}[2]{%
  \href{http://www.ams.org/mathscinet-getitem?mr=#1}{#2}
}
\providecommand{\href}[2]{#2}

\end{document}